\documentclass{amsart}
\usepackage[english]{babel}
\input xy
\xyoption{all}
\usepackage{amssymb,amsmath,amsthm}

\newcommand{\abs}[1]{\left\vert#1\right\vert}
\newcommand{\R}{\mathbb{R}}
\newcommand{\Q}{\mathbb{Q}}

\newcommand{\lie}[1]{\mathfrak{#1}}     
\newcommand{\Lie}{\mathcal{L}}          
\newcommand{\Z}{\mathbb{Z}}
\newcommand{\C}{\mathbb{C}}
\newcommand{\hook}{\lrcorner\,}
\newcommand{\LieG}[1]{\mathrm{#1}}      
\newcommand{\Spin}{\mathrm{Spin}}
\newcommand{\SU}{\mathrm{SU}}
\newcommand{\SO}{\mathrm{SO}}
\newcommand{\Gtwo}{\mathrm{G}_2}

\newcommand{\su}{\mathfrak{su}}

\newcommand{\Id}{\mathrm{Id}}
\DeclareMathOperator{\Stab}{Stab}
\DeclareMathOperator{\Int}{Int}

\newcommand{\dfn}[1]{\emph{#1}}
\newcommand{\hs}[1]{#1}
\newcommand{\ev}[1]{\mathbf{#1}}

\theoremstyle{plain}
\newtheorem{proposition}{Proposition}
\newtheorem{theorem}[proposition]{Theorem}
\newtheorem{lemma}[proposition]{Lemma}

\theoremstyle{definition}

\theoremstyle{remark}
\newtheorem*{remark}{Remark}

\newcommand{\e}{\Id}
\begin{document}
\author{Diego Conti}
\address{Dipartimento di Matematica e Applicazioni, Universit\`a di Milano  Bicocca\\  Via Cozzi 53\\ 20125 Milano\\Italy}
\email{diego.conti@unimib.it}
\title{Cohomogeneity One Einstein-Sasaki $5$-manifolds}
 \begin{abstract}  We consider hypersurfaces in Einstein-Sasaki $5$-manifolds which are tangent to the
characteristic vector field. We introduce evolution equations  that can be used to reconstruct the $5$-dimensional metric from such a hypersurface, analogous to the (nearly) hypo and  half-flat evolution equations in higher dimensions. We use these equations to classify Einstein-Sasaki $5$\nobreakdash-manifolds of cohomogeneity one.
\end{abstract}
\maketitle
\section*{Introduction}
From a Riemannian point of view, an Einstein-Sasaki manifold is a Riemannian manifold $(M,g)$ such that the
conical metric on $M\times\R^+$ is K\"ahler and Ricci-flat. In particular, this implies that $(M,g)$ is
odd-dimensional, contact and Einstein with positive scalar curvature. The Einstein-Sasaki manifolds that are simplest to describe are the \dfn{regular} ones, which arise as circle bundles over K\"ahler-Einstein manifolds.
In five dimensions, there is a classification of regular Einstein-Sasaki manifolds \cite{FriedrichKath}, in which precisely two homogeneous examples appear, namely the sphere $S^5$ and the Stiefel manifold
\begin{equation}
\label{eqn:HomogeneousSpace} V_{2,4}=\SO(4)/\SO(2)\cong S^2\times S^3\;.
\end{equation}
In fact these examples are unique (up to finite cover), as homogeneous contact manifolds are necessarily regular \cite{BoothbyWang}. Among  regular Einstein-Sasaki $5$\nobreakdash-manifolds, these two are the only ones for which the metric is known explicitly: indeed, the sphere is equipped with the standard metric, and the metric on $V_{2,4}$ has been described in \cite{Tanno}. Notice however that both $S^5$ and $S^2\times S^3$ carry other, non-regular  Einstein-Sasaki metrics \cite{BoyerGalickiKollar,BoyerGalickiNakamaye}.

Only recently other explicit examples of Einstein-Sasaki manifolds have been found, as in \cite{GauntlettMartelliSparksWaldram} the
authors constructed an infinite family $Y^{p,q}$ of Einstein-Sasaki metrics on \mbox{$S^2\times S^3$} (see also \cite{CveticEtAl:NewEinsteinSasaki} for a generalization). These metrics are non-regular, and thus they are not included in  the above-mentioned classification. The isometry group of each $Y^{p,q}$ acts with cohomogeneity one, meaning that generic orbits are hypersurfaces. In this paper we give an alternative construction of the $Y^{p,q}$, based on the language of cohomogeneity one manifolds. In fact we prove that, up to finite cover, they are the only Einstein-Sasaki $5$-manifolds on
which the group of isometries acts with cohomogeneity one.  In particular, this result settles a question raised in \cite{Gauntlett:Obstructions}, concerning the family of links $L(2,2,2,k)$, $k>0$ defined by the polynomial
\[z_1^2+z_2^2+z_3^2+z_4^k\;.\]
The homogeneous metrics mentioned earlier provide examples of Einstein-Sasaki metrics on $L(2,2,2,k)$, $k=1,2$ such that the integral lines of the characteristic vector field are the orbits of the natural action of $\LieG{U}(1)$ with weights $(k,k,k,2)$. The authors of \cite{Gauntlett:Obstructions} show that for $k>3$ no such metric exists. Since $L(2,2,2,3)$ is diffeomorphic to $S^5$, and a  metric on $L(2,2,2,3)$ of the required type would necessarily be of cohomogeneity one, our classification extends this result to $k\geq 3$.

In five dimensions, another characterization of Einstein-Sasaki manifolds is by the existence of a real Killing spinor or, more precisely, a spinor satisfying
 \begin{equation*}
 \nabla_X \psi = -\frac12 X\cdot \psi\;,
 \end{equation*}
 where $\nabla$ is the Levi-Civita connection and $\cdot$ denotes Clifford multiplication.
Such a spinor exists on every Einstein-Sasaki manifold in any dimension. The converse only fails to be true
to the extent that a complete simply-connected Riemannian manifold which admits a real Killing spinor but no
Einstein-Sasaki structure is either a sphere $S^{2n}$, a $7$-manifold with a nearly-parallel $\Gtwo$
structure  or a nearly-K\"ahler $6$-manifold~\cite{Bar}. In all of these cases the metric is Einstein with
positive scalar curvature. Thus, Einstein-Sasaki manifolds can be viewed as part of a slightly more general
class.  We shall not make much use of the spinor formalism in this paper, because once one fixes the
dimension, the spinor can be replaced with differential forms. However, the characterization in terms of
spinors establishes an analogy which is suggestive of the fact that the methods of this paper can be adapted
to the nearly-K\"ahler case, potentially leading to the construction of non-homogeneous nearly-K\"ahler
$6$-manifolds. This expectation is supported  by the fact that the corresponding exterior differential system
is involutive in  the nearly-K\"ahler  as in the Einstein-Sasaki case.  As far as the author knows, the only
known examples of nearly-K\"ahler $6$-manifolds are homogeneous (see \cite{Butruille} for a classification).
On the other hand, cohomogeneity one nearly-parallel $\Gtwo$ manifolds are classified, and none of them is
complete \cite{CleytonSwann}.

Our method to classify cohomogeneity one Einstein-Sasaki $5$-manifolds consists in considering a generic orbit, which by hypothesis
is a hypersurface. We prove that it is tangent to the characteristic vector field; this fact enables us to
define an induced global frame on the hypersurface, and write down certain equations that it must satisfy.
Resuming the analogy with higher dimensions, these equations correspond to the nearly hypo and nearly
half-flat equations of \cite{Fernandez:NearlyHypo}. In general, a hypersurface in an Einstein-Sasaki
$5$-manifold (which we take to be tangent to the characteristic vector field) determines the Einstein-Sasaki
structure locally. We express this fact in terms of evolution equations like in
\cite{Hitchin:StableForms,ContiSalamon,Fernandez:NearlyHypo}. These evolution equations are of independent
interest, as their solutions correspond to local Einstein-Sasaki metrics; the fact that a solution always
exists for real analytic initial conditions is not obvious, but it follows from the exterior differential
system being involutive. However, this is not essential to our classification: all we need to do is find all
homogeneous solutions of our nearly-hypo-like equations, and solve the evolution equations explicitly in an interval $(t_-,t_+)$, with
these solutions as initial conditions. Having done that, we prove that for suitable choices of the parameters
the resulting metric can somehow be extended to the boundary of $(t_-,t_+)$, leading to cohomogeneity one manifolds with an invariant
Einstein-Sasaki metric. Conversely, all such Einstein-Sasaki $5$-manifolds are obtained this way. This
construction is similar to the one of \cite{CleytonSwann}, where cohomogeneity one manifolds with (weak) holonomy
$\Gtwo$ were classified. However, such manifolds are  non-compact, in sharp contrast with our
case. Indeed, the metrics we obtain appear to be the first compact examples obtained by an evolution in the sense of
\cite{Hitchin:StableForms}.

\section{Invariant Einstein-Sasaki $\SU(2)$-structures }
\label{sec:SU2}
In the usual terminology, an Einstein-Sasaki structure on a $5$-manifold is a type of $\LieG{U}(2)$-structure. However, simply connected Einstein-Sasaki manifolds carry a real Killing spinor \cite{FriedrichKath},
which reduces the structure group to $\SU(2)$. The relevant representation of $\SU(2)$ can be described by the diagram
\[
\xymatrix{\SU(2)_-\ar@<-.5ex>@{^{(}->}[r]\ar@{=}[d] & \SU(2)_+\times\SU(2)_-\ar[r]^{\qquad\cong}\ar[d]^{2:1}&\Spin(5)\ar[d]^{2:1}\\
\SU(2)\ar[r]&\SO(4)\ar[r]&\SO(5)}
\]
Giving an $\SU(2)$-structure on a $5$-manifold $\ev{M}$ is the same as giving  differential forms $(\alpha,\omega_1,\omega_2,\omega_3)$,
such that locally there exists a basis of orthonormal forms $e^1,\dots,e^5$ satisfying
\begin{equation}
\label{eqn:referenceformsSU2}\left\{
\begin{aligned}\alpha&=e^5&\omega_1&=e^{12}+e^{34}\\
\omega_2&=e^{13}+e^{42}&\omega_3&=e^{14}+e^{23}
\end{aligned}\right.
\end{equation}
Here and in the sequel, we abbreviate $e^1\wedge e^2$ as $e^{12}$, and so on.
By \cite{ContiSalamon}, the Einstein-Sasaki condition can be written
\begin{align}
 \label{eqn:EinsteinSasaki}
d\alpha&=-2\,\omega_1\;, & d\omega_2&=3\alpha\wedge\omega_3\;,& d\omega_3&=-3\alpha\wedge\omega_2\;.
\end{align}
By this we mean that, up to passing to the universal cover, every Einstein-Sasaki $\LieG{U}(2)$-structure on a $5$-manifold has an $\SU(2)$-reduction satisfying \eqref{eqn:EinsteinSasaki}.

\begin{remark}
The constant $3$ appearing in Equation \eqref{eqn:EinsteinSasaki} is in some sense not essential: one could
replace it with an unspecified constant, obtaining a possible definition of an
$\alpha$\nobreakdash-Einstein-Sasaki $\SU(2)$-structure. Most of our arguments generalize to this more
general setting in a straightforward way. However, since this paper is mainly concerned with cohomogeneity
one metrics, we shall focus on the Einstein case, as the generalization does not seem to produce any new
example.
\end{remark}

We are interested in Einstein-Sasaki manifolds $\ev{M}$ of cohomogeneity one, namely those for which the principal orbits of the isometry group are hypersurfaces. Since we require $\ev{M}$ to be complete, and  Einstein-Sasaki manifolds have positive Ricci, by Myers' theorem $\ev{M}$ will be compact with finite fundamental group. We shall assume that $\ev{M}$ is simply connected, which amounts
to passing to the universal cover. In this hypothesis, we now prove some facts that play an importan r\^ole in the classification.
Recall that the \dfn{characteristic vector field} is by definition the vector field dual to $\alpha$.
\begin{lemma}\label{lemma:OrbitGo}
Let a compact Lie group $G$ act on a contact manifold $\ev{M}$ with cohomogeneity one, preserving the
contact form. Then the characteristic vector field is tangent to each principal orbit.
\end{lemma}
\begin{proof}
Let $\alpha$ be the contact form, and let $\ev{M}$ have dimension $2n+1$. On every principal orbit $Gx$ we have
\[(d\alpha)^n=d(\alpha\wedge (d\alpha)^{n-1})\;;\]
by Stokes' Theorem,
\[0=\int_{Gx} (d\alpha)^n\;.\]
On the other hand, $(d\alpha)^n$ is invariant under $G$, so it must vanish identically. Now observe that at
each point $x$, the characteristic direction is the space
\[\{X\in T_x\ev{M} \mid X\hook (d\alpha)^n=0\}\;;\]
the statement follows immediately.
\end{proof}

\begin{lemma}
\label{lemma:GPreservesSU2} Let $\ev{M}$ be a compact, simply-connected, Einstein-Sasaki $5$-manifold. Suppose
that the group of isometries of $\ev{M}$ acts with cohomogeneity one. Then the action of its identity component
$\mathcal{I}$ preserves the Einstein-Sasaki $\LieG{U}(2)$\nobreakdash-structure of $\ev{M}$, and one can define a
homomorphism
\[e^{i\gamma}\colon\mathcal{I}\to\LieG{U}(1),\quad L_g^* (\omega_2+i\omega_3)=e^{i\gamma}(g)(\omega_2+i\omega_3)\;.\]
If $\mathcal{I}$ has dimension greater than four, $\mathcal{I}'=\ker e^{i\gamma}$ is a $4$-dimensional Lie group that acts on $\ev{M}$ with cohomogeneity one,
beside preserving the Einstein-Sasaki $\SU(2)$-structure.
\end{lemma}
\begin{proof}
The action of $\mathcal{I}$ preserves both metric and orientation on $\ev{M}$, and therefore the spin structure. Thus, $\mathcal{I}$ acts on the
space of spinors $\Gamma(\Sigma)$, and the space $\mathcal{K}\subset \Gamma(\Sigma)$ of Killing spinors with
Killing constant $-1/2$ is preserved by $\mathcal{I}$.  By the assumption on the isometry group, $\ev{M}$ is not isometric to the sphere. Hence, $\mathcal{K}$ is a complex vector
space of dimension one (see \cite{FriedrichKath}), and it determines the $\LieG{U}(2)$-structure of $\ev{M}$. The action of $\mathcal{I}$ on
$\mathcal{K}\cong\C$ determines the homomorphism $e^{i\gamma}$.

Suppose $\mathcal{I}$ has dimension greater than four, and let $K$ be the stabilizer at a point of a principal orbit $M$; by dimension count, $K$ is not discrete. We claim that
\begin{equation}
\label{eqn:I0CohomogeneityOne}
\lie{k}\oplus\lie{i}'=\lie{i}\;,
\end{equation}
where $\lie{k}$, $\lie{i}$, $\lie{i}'$ are the Lie algebras of $K$, $\mathcal{I}$ and $\mathcal{I}'$ respectively.
We shall prove by contradiction that $\lie{k}$ is not contained in $\lie{i}'$, which implies \eqref{eqn:I0CohomogeneityOne} by a dimension count.

The principal orbit $M$ has trivial normal bundle, and the unit normal is invariant. Since the action of $\mathcal{I}$ preserves $\mathcal{K}$, it also preserves the contact
form, and so by Lemma~\ref{lemma:OrbitGo} the characteristic
vector field is tangent to $M$. Consequently, the unit normal is an invariant section of $\ker\alpha|_M$. The structure group $\SU(2)$ acts freely on unit vectors in $\R^4=\ker e^5$, and so $\mathcal{I}'$ preserves an $\{\e\}$-structure on $M$. Hence
$\lie{k}\subset \lie{i}'$ acts trivially on $\lie{i}/\lie{k}$, implying that $\lie{k}$ is an ideal of $\lie{i}$. Then
the identity component of $K$ is a normal subgroup of $\mathcal{I}$, and so it acts trivially on $\ev{M}$; since $\mathcal{I}$ acts effectively, this implies that $K$ is discrete, which is absurd. Thus, \eqref{eqn:I0CohomogeneityOne} holds and $\mathcal{I}'$ acts with cohomogeneity one.

By the same token, the stabilizer $\mathcal{I}'\cap K$ has trivial isotropy representation and is therefore normal in $\mathcal{I}'$. But then it is also normal in $\mathcal{I}$, which acts effectively, implying that $\mathcal{I}'\cap K$ is the trivial group and $\mathcal{I}'$ acts freely on $M$. In particular,
$\mathcal{I}'$ has dimension equal to four.
\end{proof}
\begin{remark}
Part of the argument of Lemma~\ref{lemma:GPreservesSU2} was exploited in \cite{NagyMoroianuSemmelmann} to
prove that unit Killing fields do not exist on nearly-K\"ahler $6$-manifolds other than $S^3\times S^3$.  The
$6$\nobreakdash-\hspace{0pt}dimensional case is simpler in this respect, because Killing spinors form a
\emph{real} vector space of dimension one. Thus, a Lie group acting isometrically on a nearly-K\"ahler
$6$\nobreakdash-manifold automatically preserves the nearly-K\"ahler $\SU(3)$-structure.
\end{remark}
It will follow from our classification that the isometry group of $\ev{M}$ has dimension five, but we cannot prove it directly. The best that we can do at the moment is the following:
\begin{proposition}
\label{proposition:SU2U1}
If $\ev{M}$ is a simply-connected compact Einstein-Sasaki $5$-manifold of cohomogeneity one, then the Lie group $\SU(2)\times\LieG{U}(1)$ acts on $\ev{M}$ with cohomogeneity one preserving the Einstein-Sasaki $\LieG{U}(2)$-structure, and each principal stabilizer $K$ is a finite subgroup of a maximal torus $\LieG{T}^2\subset\SU(2)\times\LieG{U}(1)$.
\end{proposition}
\begin{proof}
Since the isometry group acts with cohomogeneity one, its dimension is greater than or equal to four. If greater, by Lemma~\ref{lemma:GPreservesSU2} there is a four-dimensional subgroup that acts with cohomogeneity one.
Thus, there is a four-dimensional compact Lie group $G$ that acts with cohomogeneity one. Up to finite cover, $G$ is either a torus or $\SU(2)\times\LieG{U}(1)$. Under the projection $G\to G/K\cong M$, where $M$ is a principal orbit, the invariant one-form $\alpha$ pulls back to a non-closed left-invariant one-form on $G$. No such form exists on a torus, and so we can assume that $G$ is $\SU(2)\times\LieG{U}(1)$.

Now, the stabilizer $K$ is discrete, because $G$ acts with cohomogeneity one, and compact, because $G$ is; hence, $K$ is finite. The tangent space of a principal orbit splits as
\[TM= G\times_K(\su(2)\oplus\lie{u}(1))= G\times_K\su(2)\oplus G\times_K\lie{u}(1)\]
into two integrable distributions. Since $\ker\alpha$ is not integrable on $M$,
\[\ker\alpha|_M\neq G\times_K\su(2)\;.\]
Thus $\alpha$ has a non-zero component in $G\times_K(\su(2))^*$. On the other hand, the subgroup of $G$ that fixes a non-zero element in $\su(2)^*$ is a maximal torus $\LieG{T}^2$. Since $\alpha$ is fixed by the isotropy representation, it follows that $K\subset \LieG{T}^2$.
\end{proof}

\section{Cohomogeneity one}
\label{sec:CohomogeneityOne}
In this section we recall some standard facts and notation concerning cohomogeneity one manifolds, referring to \cite{Bredon} for the details. Like in the statement of Proposition~\ref{proposition:SU2U1}, in this section $\ev{M}$ is a compact simply-connected manifold on which a Lie group $G$ acts with cohomogeneity one.

By the general theory of cohomogeneity one manifolds, the orbit space $\ev{M}/G$ is a one-dimensional manifold, which by our topological assumptions is a closed interval. More precisely, we can fix an interval $I=[t_-,t_+]$ and a geodesic $c\colon I\to \ev{M}$  that intersects principal orbits orthogonally, so that the induced map
\begin{equation*}
I \to \ev{M}/G
\end{equation*}
is a homeomorphism. The stabilizer of $c(t)$ is fixed on the interior of $I$, and we can define
\[H_\pm=\Stab c(t_\pm),\quad K=\Stab c(t),\; t\in (t_-,t_+)\;.\]
Then $H_\pm\supset K$, and moreover there is a sphere-transitive orthogonal representation $V_\pm$ of $H_\pm$ with principal stabilizer $K$. By the simply-connectedness assumption, $V_\pm$ has dimension greater than one, so that $H_\pm/K$ is diffeomorphic to a sphere of positive dimension.

The orbits $Gc(t_\pm)\cong G/H_\pm$ are called \dfn{special}. In a neighbourhood of each special orbit, $\ev{M}$ is $G$-equivariantly diffeomorphic to an associated bundle
 \begin{equation}
 \label{eqn:AssociatedVB}
G\times_{H_\pm} D_\pm\subset G\times_{H_\pm} V_\pm\;,
 \end{equation}
where $D_\pm$ is the closed unit disk in $V_\pm$.
One can reconstruct $\ev{M}$ by glueing together two disk bundles of the form \eqref{eqn:AssociatedVB} along the boundaries, namely along the sphere bundles $G\times_{H_\pm} \partial D_\pm$. By hypothesis these boundaries consist of a single principal orbit; identifying both of them with $G/K$, the glueing map is defined by a $G$-equivariant automorphism of $G/K$. All such automorphisms have the form
\[gK\to gaK,\quad a\in N(K)\;.\]
It is customary to represent any cohomogeneity one manifold obtained this way
by the diagram \begin{equation}\label{eqn:Diagram} K\subset\{H_-,H_+\}\subset G\;.\end{equation}
Strictly speaking, this diagram does not determine $\ev{M}$ (up to $G$-equivariant diffeomorphism), because the glueing map is also involved. However in many cases, e.g. if $N(K)$ is connected, any glueing map can be extended to an equivariant diffeomorphism of at least one of the two disk bundles, and so all maps give the same manifold up to equivariant diffeomorphism (see \cite{Uchida} for details).

\smallskip
More specifically to our case, by Proposition~\ref{proposition:SU2U1} $K$ is a finite subgroup of $T^2\subset G$. Since $H_\pm/K$ is diffeomorphic to a sphere, each $H_\pm$ can be either one-dimensional or three-dimensional. However, since $\ev{M}$ is simply-connected, $H_+$ and $H_-$ cannot both be three-dimensional (see e.g. \cite{GroveWilkingZiller}).
If $H_\pm$ is three-dimensional, then its identity component is $(H_\pm)_0=\SU(2)\times\{1\}$ and $H_\pm=K\cdot (H_\pm)_0$. Now
\[S^3\cong \frac{H_\pm}{K}=\frac{(H_\pm)_0}{(H_\pm)_0\cap K}\]
and thus $K\cap (\SU(2)\times\{1\})$ is the trivial group. It follows that the tube about the special orbit $G/H_\pm$ satisfies
\begin{equation}
\label{eqn:TrivialSliceRepresentation}
G\times_{H_\pm} V_\pm= G\times_{(H_\pm)_0} V_\pm\;.
\end{equation}
Concerning the case that $H_\pm$ is one-dimensional, we remark that if one assumes that $H_\pm$ is contained in $T^2$, then $K$ is normal in $H_\pm$. Then $K$ acts trivially in the slice representation, and again \eqref{eqn:TrivialSliceRepresentation} holds. These facts will be used in the proof of Lemma~\ref{lemma:Extends}.

\section{Evolution and hypersurfaces}
\label{sec:Evolution}
We now introduce a convenient language to determine the cohomogeneity one metrics in terms of the invariant structure induced on a principal orbit. This section is motivated by Section~\ref{sec:SU2}, but otherwise independent, and should be read in the context of evolution equations in the sense of \cite{Hitchin:StableForms}.

In this section, $M$ is an oriented hypersurface in a $5$-manifold $\ev{M}$ with an $\SU(2)$-structure; notice that we make no hypotheses of invariance. Let $X$ be the unit normal to $\hs{M}$, compatible with the orientations. One can measure the amount to which $\alpha$ fails to be tangent to $\hs{M}$ by the angle
\begin{equation*}
\beta\colon\hs{M}\to[-\pi/2,\pi/2]\;,\quad \sin\beta=\iota^*(\alpha(X))\;.
\end{equation*}
The relevant case for our classification is $\beta\equiv0$ (Lemma~\ref{lemma:OrbitGo}).
\begin{remark}
If $\alpha$ is a contact form on $\ev{M}$, the angle $\beta$ cannot equal $\pm\pi/2$ on an open subset of $\hs{M}$, as
that would mean that the distribution $\ker\alpha$ is integrable.
\end{remark}
Given a hypersurface $\iota\colon\hs{M}\to\ev{M}$, one can look for deformations of $\iota$ that leave $\beta$ unchanged. Under suitable integrability assumptions on $\ev{M}$, there is a canonical  deformation obtained by the exponential map with this property. We shall say that an $\SU(2)$-structure is  \dfn{contact} if the  underlying almost contact metric structure is contact, i.e. $d\alpha=-2\omega_1$, and it is  \dfn{K-contact} if in addition the characteristic vector field is Killing. There is a well-known characterization of $K$-contact structures \cite{Blair}, which in
our language reads
\[\nabla_X\alpha=-2X\hook\omega_1\quad \forall X\in T\ev{M}\;,\]
where $\nabla$ is the Levi-Civita connection.
\begin{lemma}
\label{lemma:KContact} Let $\ev{M}$ be a K-contact 5-manifold and let $\iota\colon\hs{M}\to\ev{M}$ be an oriented, compact
 embedded hypersurface. Consider the one-parameter family of immersions
$\iota_t\colon\hs{M}\to\ev{M}$ given by
$$\iota_t(x)=\exp_{\iota(x)}(tX_x)\;,$$ where $X$ is the unit normal. Then the angle $\beta_t\colon M\to[-\pi/2,\pi/2]$ of the hypersurface $\iota_t\colon\hs{M}\to\ev{M}$ does not depend on $t$.
\end{lemma}
\begin{proof}
Using the exponential map, we can extend $X$ to a neighbourhood of $\hs{M}$, in such a way that $X$ is normal to all of
the $\iota_t(\hs{M})$ for small $t$. We must prove
 \begin{equation}
\label{eqn:alphatangent}
\iota^*_t(\alpha(X))=\iota^*(\alpha(X))\;,
\end{equation}
which holds trivially for $t=0$. By the definition of the exponential map, $\nabla_X X=0$. So,
\[\Lie_X(\alpha(X))=(\nabla_X\alpha)(X)=(-2X\hook\omega_1)(X)=0\;,\]
and \eqref{eqn:alphatangent} holds for all $t$.
\end{proof}
With notation from Lemma~\ref{lemma:KContact}, the exponential map produces an inclusion
\[\hs{M}\times(t_-,t_+)\ni (x,t)\to \iota_t(x)\in\ev{M}\;;\]
with respect to which the Riemannian metric $\ev{g}$ on $\ev{M}$ pulls back to a metric in the ``generalized cylinder form'' $dt^2+\iota_t^*\ev{g}$. Using Lemma~\ref{lemma:KContact}, we can do the same for the $\SU(2)$-structure:
\begin{proposition}\label{prop:OneParameterFamily}
Let $\iota\colon\hs{M}\to\ev{M}$ be a compact, oriented  hypersurface with angle $\beta=0$ in a $K$-contact
$5$\nobreakdash-manifold $(\ev{M},\alpha,\omega_i)$. Then there is a one-parameter family of
$\{\e\}$\nobreakdash-\hspace{0pt}structures $(\eta^0(t),\dots,\eta^3(t))$ on $\hs{M}$ satisfying
\begin{align*}
\eta^0(t)&=\iota_t^*\alpha \;,& \eta^{23}(t)&=\iota_t^*\omega_1\;,& \eta^{31}(t)&=\iota_t^*\omega_2\;,& \eta^{12}(t)&=\iota_t^*\omega_3\;,
\end{align*}
and $\ev{M}$ is locally given as the product $\hs{M}\times(t_-,t_+)$, with $\SU(2)$-structure determined by
\begin{align*}
\alpha&=\eta^0(t)\;, & \omega_1&=\eta^{23}(t) + \eta^1(t)\wedge dt\;,\\ \omega_2&=\eta^{31}(t) + \eta^2(t)\wedge dt\;,& \omega_3&=\eta^{12}(t) + \eta^3(t)\wedge dt\;.
\end{align*}
\end{proposition}
\begin{proof}
Let $x$ be a point of $\iota_t(\hs{M})$. Choose a basis $e^1,\dots,e^5$ of $T^*_x\ev{M}$ such that Equations
\ref{eqn:referenceformsSU2} hold. Using the metric on $\ev{M}$, we can write
 \[T^*_x\ev{M}=T^*_x\iota_t(\hs{M})\oplus \langle dt\rangle\;,\]
 where $dt$ represents the unit normal $1$-form compatible with the orientations.
By Lemma~\ref{lemma:KContact}, $e^5$ lies in $T^*_x\iota_t(\hs{M})$. We can act on $e^1,\dots,e^5$ by some element of $\SU(2)$ to obtain $e^4=dt$. Then
\begin{equation*}
\eta^0=e^5\;,\quad \eta^1=e^3\;,\quad \eta^2=-e^2\;,\quad \eta^3=e^1\;,\quad dt=e^4.\qedhere
\end{equation*}
\end{proof}

It is well known that Einstein-Sasaki manifolds are $K$-contact. In the hypotheses of Proposition \ref{prop:OneParameterFamily}, assume that $(\alpha,\omega_i)$ satisfy \eqref{eqn:EinsteinSasaki}. Then
it is clear that for all $t$, the following hold:
\begin{align}
\label{eqn:go}
d\eta^0&=-2\eta^{23}\;, &     d\eta^{31}&=3\eta^{012}\;, &   d\eta^{12}&=- 3\eta^{031}\;.
\end{align}
Thus, every oriented hypersurface in a $5$-dimensional Einstein-Sasaki manifold has a natural $\{\e\}$\nobreakdash-structure satisfying \eqref{eqn:go}.
 By Proposition~\ref{prop:OneParameterFamily}, a one-parameter family of $\{\e\}$\nobreakdash-structures satisfying \eqref{eqn:go} is induced on the hypersurface.
Conversely, we have the following:
\begin{proposition} \label{prop:GoEvolution}
Let $(\eta^i(t))$ be a 1-parameter family of $ \{\e\}$-structures on $\hs{M}$ such that
\eqref{eqn:go} holds for $t=t_0$. The
induced $\SU(2)$-structure on $\ev{M}=M\times (t_-,t_+)$
 is Einstein-Sasaki if and only if
 \begin{equation} \label{eqn:GoEvolution}\left\{\begin{aligned}
\partial_t \eta^0&=2\eta^1 & &&   \partial_t\eta^{23}&=-d\eta^1&\\    \partial_t\eta^{31}&=3\eta^{03}-d\eta^2 &  &&  \partial_t\eta^{12}&=-3\eta^{02}-d\eta^3
\end{aligned}\right.
 \end{equation}
In this case, \eqref{eqn:go} holds for all $t$.
\end{proposition}
\begin{proof}
Suppose Equations \ref{eqn:GoEvolution} define an Einstein-Sasaki structure on $\ev{M}$; then \eqref{eqn:go} holds for all $t$. Equations \ref{eqn:EinsteinSasaki} give
  \[dt\wedge\partial_t\eta^0+d\eta^0=d\alpha=-2\omega_1=-2\eta^{23}-2\eta^1\wedge dt\;;\]
hence $\partial_t\eta^0=2\eta^1$. Moreover
$$d\eta^{31}+dt\wedge\partial_t\eta^{31}+d\eta^2\wedge dt=d\omega_2=3\alpha\wedge\omega_3=3\eta^{012}+3\eta^{03}\wedge dt\;,$$
and therefore
 \[\partial_t\eta^{31}=3\eta^{03}-d\eta^2\;;\]
  similarly, for $\omega_3$ we get
\[\partial_t\eta^{12}=-3\eta^{02}-d\eta^3\;.\]
Since \eqref{eqn:EinsteinSasaki} implies that $\omega_1$ is closed, we also obtain
$d\eta^1=-\partial_t\eta^{23}$. So the ``only if'' part is proved.

To prove the ``if'' part, it suffices to show that (\ref{eqn:GoEvolution}) forces \eqref{eqn:go} to hold  for
all $t$;
 \eqref{eqn:EinsteinSasaki} will then follow from the calculations in the first part of the proof.
Observe first that
\[\partial_t (d\eta^0+2\eta^{23})=d\partial_t\eta^0 +2\partial_t\eta^{23}=0\;,\]
so the condition $d\eta^0=-2\eta^{23}$ holds for all $t$. Using this, we compute
\[\partial_t(d\eta^{12}+3\eta^{031})=-3d\eta^{02}-3\eta^0\wedge
d\eta^2=-3d\eta^0\wedge\eta^2=6\eta^{23}\wedge\eta^2=0\] and a similar argument shows that all of
\eqref{eqn:go} are  preserved in time.
\end{proof}
\begin{remark}
The evolution equations \eqref{eqn:GoEvolution} are not in Cauchy-Kowalewsky form. Indeed, write
$\eta^i(t)=a^i_j(x,t)\eta^j$, and $d\eta^i=c^i_{jk}\eta^{jk}$, where $c^i_{jk}$ are functions of $x$. Then
\eqref{eqn:GoEvolution} reduce to a system of $22$  equations in $16$ unknowns; so the system is
overdetermined. However, it turns out that a solution always exists, at least in the real analytic case, as
will be proved elsewhere.
\end{remark}
\begin{remark}
Given an $\{\e\}$\nobreakdash-structure satisfying \eqref{eqn:go}, we can define another $\{\e\}$\nobreakdash-\hspace{0pt}structure satisfying \eqref{eqn:go} by
\begin{equation}
\label{eqn:signchange}
(\eta^0,\eta^1,\eta^2,\eta^3)\to (\eta^0,-\eta^1,-\eta^2,-\eta^3)\;;
\end{equation}
this has the effect of reversing the orientation, and reflects the ambiguity in the reduction from $\Z_2$ to the trivial group. In terms of the evolution equations, this change in orientation corresponds to
changing the sign of $t$.
\end{remark}
As an example, consider the Lie group $\SU(2)\times\LieG{U}(1)$. Fix a basis $e^1,e^2,e^3,e^4$  of left-invariant $1$-forms  such that
\begin{align}
\label{eqn:BasisForSU2U1}
 de^1&=-e^{23}\;,   &  de^2&=-e^{31}\;,   &  de^3&=-e^{12} \;,  &   de^4&=0\;.
\end{align}
One can then define a left-invariant $\{\e\}$\nobreakdash-structure satisfying \eqref{eqn:go} by
\begin{align*}
\eta^0&=\frac13e^1+e^4\;, &    \eta^1&=e^4\;, &  \eta^2&=\frac1{\sqrt6}\,e^2\;,   &
\eta^3&=\frac1{\sqrt6}\,e^3\;.
\end{align*}
The solution of the evolution equations is:
\begin{align}
\label{eqn:DegenerateSolution} \eta^0&=\frac13e^1+\cos\,\epsilon t\,e^4\;, &
\eta^1&=-\frac{\epsilon}2\sin\,\epsilon t\,e^4\;, & \eta^2&=\frac1{\epsilon}\,e^2\;, &
\eta^3&=\frac1{\epsilon}\,e^3\;,
\end{align}
where we have set $\epsilon=\sqrt6$. The resulting Einstein-Sasaki structure  can be extended to a compact
$5$-manifold, which can be realized as a circle bundle over the K\"ahler-Einstein manifold $S^2\times S^2$.
Indeed, this is the  homogeneous Kobayashi-Tanno metric (see \cite{Tanno}) on the  space
\eqref{eqn:HomogeneousSpace}.
   To see
this, let $S^3$ be another copy of $\SU(2)$, with global invariant forms $\tilde e^i$ satisfying relations analogous to
\eqref{eqn:BasisForSU2U1}. We can introduce coordinates $\theta,\psi,\phi$ on $S^3$ such that
\begin{equation}
\label{eqn:Euler}
\left\{
\begin{gathered}
\tilde e^1=-d\psi+\cos\theta d\phi\quad \tilde e^2=-\sin\psi\sin\theta d\phi+\cos\psi d\theta\\ \tilde e^3=\cos\psi\sin\theta d\phi+\sin\psi d\theta
\end{gathered}\right.
\end{equation}
If we set $t=\frac1{\epsilon}\theta$ and identify $e^4$ with $\frac13 d\phi$, then \eqref{eqn:DegenerateSolution}
restricted to $\SU(2)\times \{\psi=0\}$ becomes
\begin{align}
\label{eqn:HomogeneousMetric}
\eta^0&=\frac13 (e^1 +\tilde e^1)\;,& \eta^1&=-\frac1{\epsilon}\tilde e^3\;,&
\eta^2&=\frac1{\epsilon}e^2\;,&\eta^3&=\frac1{\epsilon}e^3\;,& dt&=\frac1{\epsilon}\tilde e^2\;.
\end{align}
On the other hand,  if viewed as forms on $\SU(2)\times S^3$, the forms \eqref{eqn:HomogeneousMetric} annihilate a vector field which we may denote by $e_1-\tilde e_1$, and their Lie derivative with respect to $e_1-\tilde e_1$ is zero; so, they pass onto the quotient. Thus the forms $(\alpha,\omega_i)$ defined as in Proposition~\ref{prop:OneParameterFamily} are invariant under the left action of $\SU(2)\times\SU(2)$, and they define a homogeneous
Einstein-Sasaki structure on a homogeneous space equivalent to \eqref{eqn:HomogeneousSpace}. Notice that
this is not a symmetric space.

\section{Solutions  on \protect{$\SU(2)\times \LieG{U}(1)$}}
\label{sec:Solutions} In this section we classify solutions of
\eqref{eqn:go} on \mbox{$G=\SU(2)\times \LieG{U}(1)$} which are invariant in a certain sense. In order to determine the relevant notion of invariance, let us go back to the hypotheses of Proposition~\ref{proposition:SU2U1}. Thus, $G$ acts transitively on a hypersurface $M\subset\ev{M}$; let $(\tilde\eta^i)$ be the $\{\e\}$\nobreakdash-structure induced on $M$ by the invariant Einstein-Sasaki structure. Then $(\tilde\eta^i)$ pulls back to an $\{\e\}$\nobreakdash-structure on $G$ which we also denote by
$(\tilde\eta^i)$. By Lemma~\ref{lemma:GPreservesSU2}, $L_g^*$ acts on $\tilde\eta^2+i\tilde\eta^3$ as multiplication by $e^{i\gamma}$.
So if we set
\[\tilde\eta^2+i\tilde\eta^3=e^{i\gamma}(\eta^2+i\eta^3)\;,\]
and $\eta^0=\tilde\eta^0$, $\eta^1=\tilde\eta^1$, then the $\eta^i$ are left-invariant forms on $G$.

The Lie algebra of $G$ is \mbox{$\lie{g}=\su(2)\oplus\lie{u}(1)$}, whose dual will be represented by the basis $e^1,\dots,e^4$  satisfying \eqref{eqn:BasisForSU2U1}. Since $\lie{g}$ has only one $3$-dimensional subalgebra, the kernel of $e^{i\gamma}$ contains $\SU(2)\times\{1\}$. We can therefore write $d\gamma= m e^4$ for some integer $m$. In terms of the invariant basis $(\eta^i)$,  Equations~\ref{eqn:go} become:
\begin{align}
\label{eqn:go2}
d\eta^0&=-2\eta^{23}, &     d\eta^{31}&=3\eta^{012}+me^4\wedge\eta^{12}, &   d\eta^{12}&=- 3\eta^{031}-me^4\wedge\eta^{31}.
\end{align}
\sloppy
Observe that  the space of left-invariant solutions of \eqref{eqn:go} or \eqref{eqn:go2}  is closed under right translation. Left-invariant $\{\e\}$\nobreakdash-structures on $G$ can be identified with $\{\e\}$\nobreakdash-structures on the Lie algebra $\lie{g}$, and right translation on $G$ correponds to the adjoint action of $G$ on
$\lie{g}$. Observe also that replacing $e^{i\gamma}$ with $e^{i(\gamma+A)}$ does not affect \eqref{eqn:go}. Hence solutions of \eqref{eqn:go2} are closed under the transformation
\begin{equation}
\label{eqn:U1Action}
(\eta^0,\eta^1,\eta^2,\eta^3)\to \left(\eta^0,\eta^1,c\eta^2+s\eta^3,-s\eta^2+c\eta^3\right)
\end{equation}
where  $c$ and $s$ are real constants with $c^2+s^2=1$.

\fussy The action of $\Int{\lie{g}}$ and  transformations of type \eqref{eqn:U1Action}, \eqref{eqn:signchange} commute;  they
generate an equivalence relation on the space of solutions.
\begin{proposition}\label{prop:GoClassification}
Every invariant solution of \eqref{eqn:go2}  on
$G$ can be written up to equivalence as either
\begin{align}
\label{eqn:GoGivingNothing} \eta^0&=2hk\,e^1-\frac m3\,e^4\;, &    \eta^1&=a\,e^1\;, &  \eta^2&=h\,e^2\;,   &
\eta^3&=c\,e^2+k\,e^3\;,
\end{align}
where $h>0$, $k>0$, $a>0$ and $c$ is an arbitrary constant, or
\begin{align}
\label{eqn:GoGivingYpq}
\eta^0&=2h^2e^1+\mu \,e^4\;, &    \eta^1&=a_1e^1+a_4e^4\;, &  \eta^2&=h\,e^2\;,   & \eta^3&=h\,e^3\;,
\end{align}
where $h>0$, $a_4\neq 0$ and $3a_1\mu=6h^2 a_4-a_4-a_1m$.
\end{proposition}
\begin{proof}
Up to an inner automorphism, we can assume that $\eta^0$ is in $\langle e^1,e^4\rangle$. Using
\eqref{eqn:go2}, we deduce that $\eta^{23}=-\frac12d\eta^0$ is a multiple of $e^{23}$, meaning that
 \[\langle \eta^2,\eta^3\rangle = \langle e^2,e^3\rangle\;.\]
Up to a transformation of type \eqref{eqn:U1Action}, we can assume  $\eta^2=h\,e^2$ for some positive constant $h$.

Now suppose that $\eta^{31}$ is closed, and therefore $\eta^1\in\langle e^1,e^2,e^3\rangle$. Then \eqref{eqn:go2} gives
\[3\eta^{012}+m\,e^4\wedge\eta^{12}=0\;.\]
Hence $\eta^0= -\frac m3 e^4 + q e^1$ and $\eta^1$ lies in $\langle e^1,e^2\rangle$. The same argument applied to
$3\eta^{031}+m\,e^4\wedge\eta^{31}$ shows that $\eta^1$ is a multiple of $e^1$, and so \eqref{eqn:GoGivingNothing} holds.
Up to changing the signs of both $e^1$ and $e^3$, we can assume $a>0$. Finally, we can fix the overall orientation to obtain $k>0$.

Suppose now that $\eta^{31}$ is not closed, so that
 \[\eta^1=\beta+a_4e^4\;,\quad a_4\neq 0,\; \beta\in\langle e^1,e^2,e^3\rangle\;.\]
By \eqref{eqn:go2}, $d\eta^{31}\wedge e^2$ is
zero, and therefore {$d\eta^3\wedge e^{42}=0$}, i.e. $\eta^3=k e^3$ for some constant $k$. Since we are allowed to
change the signs of both $e^1$ and $e^3$, we can assume that $k$ is positive.
On the other hand every exact $2$-form, as is $d\eta^{31}$, gives zero on wedging with $e^4$; so, by \eqref{eqn:go2} $\beta$ is in $\langle
e^1,e^2\rangle$. The same argument applied to the last of \eqref{eqn:go2}
shows that $\beta=a_1e^1$ for some constant $a_1$.

Now write $\eta^0=\mu_1 e^1+\mu_4 e^4$, where $\mu_1\neq 0$ because $\eta^0$ is not closed. Then \eqref{eqn:go2} can be
rewritten as
 \begin{align*}
hk&=\frac12 \mu_1,& -a_4k&= h(3a_1\mu_4-3\mu_1a_4+a_1m),& -a_4h&=k(3a_1\mu_4-3\mu_1a_4+a_1m);
\end{align*}
the solution is
 \[k=h\;,\quad \mu_1=2h^2\;, \quad 3a_1\mu_4=6h^2a_4-a_4-a_1m\;.\qedhere\]
\end{proof}

\section{The local metrics}
\label{sec:Explicit}
In order to give a \emph{local} classification, i.e. classify invariant Einstein-Sasaki structures on non-compact cohomogeneity one manifolds $\SU(2)\times\LieG{U}(1)\times(t_-,t_+)$, it is now sufficient to solve the evolution equations \eqref{eqn:GoEvolution} using the solutions of Section~\ref{sec:Solutions} as initial data. Observe that taking into account the $e^{i\gamma}$ rotation of Lemma~\ref{lemma:GPreservesSU2}, we must replace the bottom row of \eqref{eqn:GoEvolution} with
 \begin{align*}
 \partial_t\eta^{31}&=3\eta^{03}-d\eta^2+d\gamma\wedge \eta^3\;, &  &&  \partial_t\eta^{12}&=-3\eta^{02}-d\gamma\wedge \eta^2-d\eta^3 \;.
\end{align*}
Notice that $e^{i\gamma}$ does not depend on $t$.

We shall distinguish among three cases.

\smallskip
(\emph{i}) The first case is given by \eqref{eqn:GoGivingYpq} when $a_1=0$. Then $6h^2-1=0$, and
\begin{align*}
\eta^0&=\frac13e^1+\mu e^4\;, &    \eta^1&=a_4e^4\;, &  \eta^2&=\frac{1}{\sqrt6}\,e^2\;,   &
\eta^3&=\frac{1}{\sqrt6}\,e^3\;.
\end{align*}
This family is closed under evolution; explicitly, setting $\epsilon=\sqrt6$ for short, the ``rotated evolution equations'' are  solved by
\begin{align*}
\eta^0&=\frac13e^1+\left(k\cos \epsilon t-\frac m3\right) e^4\;, &    \eta^1&=-\frac {k\epsilon}2\sin\epsilon t\, e^4\;, &  \eta^2&=\frac{1}{\epsilon}\,e^2\;,   &
\eta^3&=\frac{1}{\epsilon}\,e^3\;.
\end{align*}
It is not hard to check that, regardless of $k$ and $m$, the resulting metric is the Kobayashi-Tanno metric on
$V_{2,4}$ described in Section \ref{sec:Evolution}.

\smallskip
(\emph{ii}) The second case is given by \eqref{eqn:GoGivingYpq} when $a_1\neq 0$. We can set $a_4=Ca_1$, so that writing $a$ for $a_1$:
\begin{align*}
\eta^0&=
2h^2e^1+\left(2C h^2-\frac{C+m}3\right)e^4 &    \eta^1&=a(e^1+C e^4) &  \eta^2&=h\,e^2   &
\eta^3&=h\,e^3
\end{align*}
This family is  closed under evolution for every choice of $C$. Indeed, the evolution equations read
  \begin{equation}
 \label{eqn:EvolutionOfSU2U1}\frac{d}{dt} h^2=a\;,\quad \frac{d}{dt}(ah)=h-6h^3\;.
 \end{equation}
which we can rewrite as
\begin{equation*}
(2h'h^2)'=-6h^3+h
\end{equation*}
Fix initial conditions $h(t_0)=h_0$, $a(t_0)=a_0$; we assume that $(\eta^i(t))$ is a well-defined
$\{\e\}$\nobreakdash-structure at $t=t_0$, and therefore $h_0\neq 0\neq a_0$. Up to changing the orientation
of the $4$-manifold --- or equivalently, changing the sign of $t$ in \eqref{eqn:EvolutionOfSU2U1} --- we can
assume that $a_0h_0$ is positive. By continuity and non-degeneracy,  $a(t)h(t)>0$ in $(t_-,t_+)$. Thus, the
first of~\eqref{eqn:EvolutionOfSU2U1} tells us that $dh/dt>0$ in $(t_-,t_+)$. In particular, we can use $h$
as the new variable; writing $x=\frac{dh}{dt}$, \eqref{eqn:EvolutionOfSU2U1} gives
\[4x^2+2hx \frac{dx}{dh}=1-6h^2\;.\]
We can solve this differential equation explicitly  about $h_0$; the solution is
\[x(h)=\frac1{2h^2}\sqrt{A+h^4-4h^6}\;,\]
 where the sign of the square root has been chosen consistently with $h_0>0$, $a_0>0$, and $A$ is a constant determined by the initial data. More precisely, by the first of~\eqref{eqn:EvolutionOfSU2U1}
\[A=4{h_0}^6-{h_0}^4+(a_0h_0)^2\;.\]
Since the function $h\to 4h^6-h^4$ has $-\frac1{108}$ as its minimum and $a_0h_0>0$, it follows that
\begin{equation*}
A>-\frac1{108}\;.
\end{equation*}
By the above argument, $h(t)$ is a solution of
\begin{equation}
\label{eqn:ChangeOfVariable} h'(t)=\frac1{2h^2}\sqrt{A+h^4-4h^6}\;;
\end{equation}
conversely, it is clear that every solution of \eqref{eqn:ChangeOfVariable} gives rise to a solution of
\eqref{eqn:EvolutionOfSU2U1}. However we do not need to solve \eqref{eqn:ChangeOfVariable}  explicitly: it is
sufficient to set
\begin{equation}
\label{eqn:ExplicitMetric}
\left\{\begin{gathered} \eta^0=2h^2(e^1+Ce^4)-\frac{C+m}3 e^4\;, \quad
\eta^1=\frac1h\sqrt{A+h^4-4h^6}(e^1+C\,e^4)\;, \\ \eta^2=h\,e^2\;,\quad \eta^3=h\,e^3\;, \quad dt =
\frac{2h^2}{\sqrt{A+h^4-4h^6}}dh\;.
\end{gathered}
\right. \end{equation}
Whether this metric extends to a compact Einstein-Sasaki $5$-manifold depends only on the initial conditions, in a way that we shall determine in Section~\ref{sec:Compactifications}.

\smallskip
(\emph{iii})
The structure \eqref{eqn:GoGivingNothing} can be rewritten as
\begin{align*}
\eta^0&=2(hk-bc)e^1-\frac m3\,e^4\;, &    \eta^1&=a\,e^1\;, &  \eta^2&=h\,e^2+b\,e^3\;,   & \eta^3&=c\,e^2+k\,e^3\;,
\end{align*}
where we can assume that $a$ and $hk-bc$ are positive. We shall also assume that at least one of $h-k$, $b$ and $c$ is non-zero, since otherwise we are reduced to case (\emph{ii}).
This family is closed under evolution: indeed, the evolution equations give
\begin{equation}
\label{eqn:EvolutionGivingNothing}
\left\{
\begin{aligned}
\partial_t(ac)&=-6c(hk-bc)-b &\quad\partial_t(ak)&=-6k(hk-bc)+h \\
\partial_t(ab)&=-6b(hk-bc)-c &\quad\partial_t(ah)&=-6h(hk-bc)+k \\
\partial_t(hk-bc)&=a
\end{aligned}\right.
\end{equation}
Any solution of these equations defines a local Einstein-Sasaki metric, which cannot however be extended to a complete metric, as we
shall prove without solving \eqref{eqn:EvolutionGivingNothing} explicitly.

\section{Extending to special orbits}
\label{sec:Extending}
In Section~\ref{sec:Explicit} we have constructed invariant Einstein-Sasaki metrics on $G\times\!(t_-,t_+)$. The problem remains of determining whether one can extend the resulting Einstein-Sasaki structure to a compact cohomogeneity one manifold with diagram \eqref{eqn:Diagram}. In this section we study this problem in slightly greater generality.

Given a differential form (or, more generally, a tensor) defined away from the special orbits, we want to
determine conditions for it to extend to all of $\ev{M}$. By Section~\ref{sec:CohomogeneityOne}, this problem reduces to
extending a form on $G\times_H (V\setminus\{0\})$ to the zero section. Fix a
left-invariant connection on $G$ as a principal bundle over $G/H$, i.e. a $H$\nobreakdash-invariant one-form
$\omega\colon\lie{g}\to\lie{h}$ extending the identity on $\lie{h}$. Then the tangent bundle of $G\times_H V$ can
be identified with
\[(G\times V)\times_H \left(\frac{\lie{g}}{\lie{h}}\oplus V\right)\;,\]
where the principal $H$ action on $G\times V$ is given by $(R_h,\rho(h^{-1}))$. Explicitly, the
identification is induced by the $H$\nobreakdash-equivariant map
\begin{equation}
\label{eqn:TangentBundle}
 T(G\times V)\ni(g,v;A,v') \to \left(A_{\lie{g}/\lie{h}}, v'-\omega(A)\cdot v\right)\in \frac{\lie{g}}{\lie{h}}\oplus V\;,
\end{equation}
where we have used the trivialization of $T(G\times V)$ given by left translation.

A differential form on $G\times_H V$  can be viewed as a $H$\nobreakdash-equivariant map
\begin{equation*}
\tau \colon G\times (V\setminus\{0\}) \to \Lambda^*\left(\frac{\lie{g}}{\lie{h}}\oplus V\right)\:;
\end{equation*}
the form is invariant if the map $\tau$ is also $G$-invariant, i.e. $\tau(gh,v)=\tau(h,v)$. So invariant
forms are $H$\nobreakdash-equivariant maps
\begin{equation*}
\tau \colon  V\setminus\{0\} \to \Lambda^*\left(\frac{\lie{g}}{\lie{h}}\oplus
V\right)=\Lambda^*\frac{\lie{g}}{\lie{h}}\otimes\Lambda^*V\;.
\end{equation*}
 Now decompose the target space of
$\tau$ into irreducible $H$\nobreakdash-modules. It is clear that $\tau$ extends smoothly across $0$ if and only if its
components extend. Thus, we only need to determine when an equivariant map
\[\tau\colon V\setminus\{0\}\to W\]
extends, where $W$ is an irreducible $H$\nobreakdash-module.

For the rest of this section we shall assume that $H\cong\LieG{U}(1)$ and $K=\Z_\sigma$ is the cyclic subgroup of order $\sigma$; this is essentially the only situation we will need to consider in this paper. The representations of $H$ are modeled on $\C$, on which we shall use both polar and Euclidean coordinates, writing the generic element as $x+iy=re^{i\theta}$. More precisely, the non-trivial irreducible representations of $H$ are
two-dimensional real vector spaces $V_n\cong\C$, on which $e^{i\theta}$ acts as
multiplication by $e^{in\theta}$.  We have
\begin{equation*}
V_n\otimes V_m=V_{n-m}\oplus V_{n+m}\;.
\end{equation*}
The $\LieG{U}(1)$-space $V_\sigma\cong \C$ has a canonical section $[0,+\infty)\subset\C$, and any equivariant map from
$V_\sigma$ to $V_n$ is determined by its restriction to this section.
\begin{proposition}[Kazdan-Warner]
\label{prop:ExtensionCriterion} A non-zero smooth map
\[\tau\colon [0,+\infty)\to V_n\]
extends to a smooth equivariant map defined on all of $V_\sigma$  if and only if $\sigma$ divides $n$ and the $k$-th derivative of $\tau$
at zero vanishes for
\[k=0,1,\dots,\abs{\frac n\sigma}-1,\abs{\frac n\sigma}+1,\abs{\frac n\sigma}+3\dots\]
\end{proposition}
\begin{proof}\
Set $m=n/\sigma$. The equivariant extension of $\tau$ to $\R\subset V_\sigma$ must satisfy
\[\tau(-r)=(-1)^m\tau(r)\;,\] implying that $m$ is integer. We need this extension to be smooth, which is
equivalent to requiring that the odd derivatives at zero vanish if $m$ is even, and the even derivatives
vanish if $m$ is odd. The equivariant extension of $\tau$ to $V_\sigma$  is given in polar coordinates by
\[f(r,\theta)=\tau(r e^{i\theta})=e^{im\theta}\tau(r)\;.\]
Set $\tilde f(r\cos\theta,r\sin\theta)=f(r,\theta)$. If $\tilde f$ is smooth at the origin, then
\[r^ke^{im\theta}\frac{d^k\tau}{dr^k}(0)=r^k\frac{\partial^k}{\partial r^k} f(0,\theta)=r^k\sum_{h=0}^k\cos^h\theta\sin^{k-h}\theta\frac{\partial \tilde f}{\partial x^h\partial y^{k-h}}(0,0)\]
is a homogeneous polynomial of degree $k$ in $r\cos\theta$, $r\sin\theta$. This condition is also sufficient for $\tilde f$ to be smooth
\cite{KazdanWarner}. On the other hand, it is easy to check that $r^ke^{im\theta}$ is not a homogeneous
polynomial of degree $k$ unless $k-\abs{m}$ is even and non-negative.
\end{proof}

We now illustrate the method in the case of $G=\SU(2)\times\LieG{U}(1)$. On the Lie algebra of $G$, fix a basis $e_1,\dots,e_4$ such that the dual basis satisfies \eqref{eqn:BasisForSU2U1}. If we identify $\SU(2)$ with unit quaternions and its Lie algebra with imaginary quaternions, our choice can be expressed as
\[e_1=i/2,\quad e_2=j/2,\quad e_3=k/2.\]
In particular, we see from this that $e_1$ has period $4\pi$.
Let $H$ be the connected subgroup of $G$ with Lie algebra spanned by
\[\xi=p e_1+q e_4\;,\]
where $p/2$ and $q$ are coprime integers. Then $\xi$ has period $2\pi$, and we can identify $H$ with $\LieG{U}(1)$. The adjoint action of $H$ on $\lie{g}^*$ is trivial on $\langle e^1,e^4\rangle$, and so its action is determined by
\[e^{i\theta} \cdot e^2=\cos p\theta\, e^2 + \sin p\theta\, e^3\;,\quad e^{i\theta}\cdot e^3=-\sin p\theta\, e^2 + \cos p\theta\, e^3\;.\]
So, for instance, the space of $G$-invariant $2$-forms is identified with the space of $\LieG{U}(1)$-invariant maps from $V_\sigma$ to
\begin{equation*}
\Lambda^2(V_p\oplus V_\sigma\oplus \R)\cong V_0\oplus V_{p+\sigma}\oplus V_{p-\sigma}\oplus V_p\oplus
V_\sigma\;.
\end{equation*}
More explicitly, fix the connection form \[\frac1{p^2+q^2} (pe^1+qe^4)\;,\]  where we have used $\xi$ to
identify $\lie{h}$ with $\R$. This choice determines a horizontal space $\lie{h}^\perp=\langle
qe_1-pe_4,e_2,e_3\rangle$, and realizes the identification map \eqref{eqn:TangentBundle} as
\begin{align*}
(\e,r;A,v')&\to (A,v')\;,\quad A\in\lie{h}^\perp\\
(\e,r;pe_1+qe_4,v')&\to (0,v'-\sigma r\partial/{\partial y})
\end{align*}
 By duality, we see that:
 \begin{itemize}
 \item the $1$-forms $e^2$, $e^3$ and $qe^1-pe^4$ are mapped to the corresponding elements of $(\lie{g}/\lie{h})^*$;
 \item the connection form is mapped to $-\frac1{\sigma r}dy$.
 \end{itemize}
This is all one needs to know in order to translate a differential form (or a tensor) into an equivariant map $\tau$, and apply the criterion of Proposition~\ref{prop:ExtensionCriterion}.

\section{From local to global}
We can now apply the criteria from Section~\ref{sec:Extending} to determine conditions for the local metrics of Section~\ref{sec:Explicit} to extend smoothly across the special orbits.
\label{sec:Compactifications}
Consider the one-parameter family of $\{\e\}$-structures
\begin{equation}
\label{eqn:GoGeneric}\left\{
\begin{gathered}
\eta^0=2\Delta (e^1+Ce^4)-\frac{C+m}3e^4 \;, \quad   \eta^1=\Delta'(e^1+Ce^4)\;, \quad \\ \eta^2=h\,e^2+b\,e^3\;,\quad    \eta^3=c\,e^2+k\,e^3\;,
\end{gathered}\right.
\end{equation}
where we have set
\[\Delta=hk-bc,\]
and we assume the constant $C+m$ and the functions $\Delta(t)$, $\Delta'(t)$ to be non-zero in the interval $(t_-,t_+)$.
In particular, solutions of the evolution equations \eqref{eqn:ExplicitMetric} or \eqref{eqn:EvolutionGivingNothing} are of this form.
Recall from Proposition~\ref{proposition:SU2U1} that the principal stabilizer $K$ is contained in a torus $T^2$; by the expression of \eqref{eqn:GoGeneric}, it follows that the Lie algebra of $T^2$ is spanned by $e_1$ and $e_4$.

We introduce the variables
\[r_\pm = \pm(t_\pm-t)\;,\]
so that the special orbit $G/H_\pm$ corresponds to $r_\pm=0$. For brevity's sake, we shall drop the subscript
sign and simply write $H$, $r$ rather than $H_\pm$, $r_\pm$. We shall refer to the point $r=0$ as the origin,
and say that a function $f(r)$ on $[0,\epsilon)$, $\epsilon>0$ is odd (resp. even) at the origin if it
extends to a smooth odd (resp. even) function on $(-\epsilon,\epsilon)$. We can now determine conditions on
the functions $h$, $k$, $b$, $c$, and $\Delta$ in order that the metric extend smoothly across each special
orbit.
\begin{lemma}
\label{lemma:Extends}
The one-parameter family of $\{\e\}$-structures \eqref{eqn:GoGeneric} defines an invariant $\LieG{U}(2)$-structure on the cohomogeneity one manifold $G\times (t_-,t_+)$. If this structure extends across the special orbit $G/H$, then
\begin{itemize}
\item $H_0=\SU(2)\times\{1\}$, $\frac1{r^4}(hb+ck)$ is even, and
\[\Delta/r^2,\quad \frac1{r^2}(h^2+c^2),\quad \frac1{r^2}(k^2+b^2)\]
are even functions, all taking the value $\frac14$ at the origin; or
\item $H_0\cong\LieG{U}(1)$, with Lie algebra $\lie{h}=\langle pe_1+qe_4\rangle$, where $p/2$ and $q$ are coprime integers satisfying
$p+qC>0$;
\[\Delta\;,\quad h^2+c^2+b^2+k^2\]
are even functions satisfying
\begin{equation*}
\abs{\Delta''(0)}=\frac\sigma{p+qC},\quad \Delta(0)=\frac{q(C+m)}{6(p+qC)}\neq 0\;,
\end{equation*}
where $\sigma$ is the order of $K\cap H_0$,
and $hk+bc$, $h^2+c^2-b^2-k^2$ are smooth functions that, if $p\neq 0$, vanish at the origin.
\end{itemize}
The first condition is also sufficient, whereas the second one becomes sufficient if $hk+bc$ and $h^2+c^2-b^2-k^2$ are identically zero.
\end{lemma}
\begin{proof}
We think of a $\LieG{U}(2)$-structure as defined by a $1$-form $\alpha$, a $2$-form $\omega_1$ and a
Riemannian metric $g$; the $\LieG{U}(2)$-structure extends  across a special orbit if and only if these
tensors extend smoothly, while remaining non-degenerate. Assume first that $H$ has dimension three. By the final remarks of Section~\ref{sec:CohomogeneityOne}, the identity component of $H$ is the subgroup $\SU(2)\times\{1\}$, and the normal bundle at $G/H$ is given by
\eqref{eqn:TrivialSliceRepresentation}.
Thus, $V$ is the $4$-dimensional irreducible representation of $\SU(2)$ and
\[G\times_{H} V\cong V\times\LieG{U}(1)\;.\]
In particular we can identify $G/K\times(t_-,t_+)$ with an open set in $V\times\LieG{U}(1)$.
Explicitly, we define $V$ to be the field of quaternions, identifying $\SU(2)$ with the unit sphere in $V$, acting
on the left.

We write the generic element of $V$ as $x^0+x^1i+x^2j+x^3k$, and identify $r$ with the radial
coordinate
\[r=\sqrt{(x^0)^2+(x^1)^2+(x^2)^2+(x^3)^2}\;.\] \fussy Then
\begin{align*}
e^1&=\frac2{r^2}\left(x^0dx^1-x^1dx^0-x^2dx^3+x^3dx^2\right)\\
e^2&=\frac2{r^2}\left(x^0dx^2+x^1dx^3-x^2dx^0-x^3dx^1\right)\\
e^3&=\frac2{r^2}\left(x^0dx^3-x^1dx^2+x^2dx^1-x^3dx^0\right)\\
\mp dt&=\frac1{r}\left(x^0dx^0+x^1dx^1+x^2dx^2+x^3dx^3\right)
\end{align*}
where the sign of $\mp dt$ depends on the special orbit $G/H_\pm$ we consider, and the coefficient $2$ is a consequence of our choice \eqref{eqn:BasisForSU2U1}.
Rewriting in these terms the $2$-form
\[\omega_1=\Delta' (e^1+Ce^4)\wedge dt +\Delta e^{23}\;,\]
we see that $\omega_1$ is smooth and non-degenerate at the origin only if $\Delta'\sim\mp\frac{2\Delta}r$, so that
\[\omega_1\sim \frac{2\Delta}{r^2}\left(2dx^2\wedge dx^3-2dx^0\wedge dx^1+Ce^4\wedge \left(x^0dx^0+x^1dx^1+x^2dx^2+x^3dx^3\right)\right)\;.\]
Hence,  $\Delta/r^2$ is even and non-zero at the origin. In particular $\alpha$ is smooth, and it does not vanish at $G/H$ since $C+m\neq 0$. Moreover \mbox{$\Delta(0)=0=\Delta'(0)$}.
The rest of the statement now follows from the fact that the metric tensor $dt\otimes dt+\sum_i \eta^i\otimes\eta^i$  has the form
\begin{multline*}
\left(4\Delta^2+(\Delta')^2\right) e^1\otimes e^1+\left(\left(2C\Delta-\frac {C+m}3\right)^2+(C\Delta')^2\right) e^4\otimes e^4+\\ +\left(C(\Delta')^2+2\Delta\left(2C\Delta-\frac {C+m}3\right)\right)e^1\odot e^4
  +(h^2+c^2)e^2\otimes e^2 +\\+ (k^2+b^2)e^3 \otimes e^3 + (hb+ck)e^2\odot e^3+dt\otimes dt
\end{multline*}
and thus smoothness requires
\[4\Delta^2+(\Delta')^2\sim h^2+c^2\sim k^2+b^2 \sim \frac14 r^2\;.\]

Now suppose $H$ is one-dimensional, so that $H_0\cong\LieG{U}(1)$. Then $\Delta$ is non-zero at the origin.
Indeed, assuming otherwise, we shall prove that $\alpha$ is not smooth. In general, we are not allowed to replace $H$ with $H_0$, since a smooth invariant form on $G\times_{H_0} V$ does not necessarily correspond to a smooth form on $G\times_{H} V$. However, the converse always holds.
Thus, we can apply the language of Section~\ref{sec:Extending} in order to prove that $\alpha$ is not smooth. Let the Lie algebra of $H_0$ be spanned by
\[\xi=pe_1+p'e_2+p''e_3+qe_4\;.\]
 The
connection form corresponding to $H_0$ has the form \[\frac{pe^1+p'e^2+p''e^3+qe^4}{p^2+(p')^2+(p'')^2+q^2}\;.\] Since $\alpha(\xi)=2p\Delta +
q(2C\Delta - \frac{C+m}3)$,  we can write
\[\alpha \equiv -\frac1{\sigma r}\left(2p\Delta + q\left(2C\Delta - \frac{C+m}3\right)\right)dy \pmod{\left(\lie{g}/\lie{h}\right)^*}\;,\]
So if $\Delta(0)$ were zero, then  $\alpha$ could only be smooth if $C+m=0$, which is absurd.

Having shown that $\Delta(0)\neq 0$, we can conclude that the isotropy representation of $H$ fixes $e^1$, and so $H$ is contained in $T^2$.
Recall from Section~\ref{sec:CohomogeneityOne} that this condition ensures that \eqref{eqn:TrivialSliceRepresentation} holds; so, it is now safe to replace $H$ with $H_0$.

Let the Lie algebra of $H$ be generated
by $\xi=pe_1+qe_4$, where $\frac p2$ and $q$ are coprime integers. Write
\[\alpha=-\frac1{\sigma r}\left(2p\Delta+q\left(2C\Delta-\frac{C+m}3\right)\right)dy+\left(2q\Delta -p\left(2C\Delta-\frac{C+m}3\right)\right)\frac{qe^1-pe^4}{p^2+q^2}\]
By Proposition~\ref{prop:ExtensionCriterion}, $\alpha$ is smooth  if and only if
\[2(p+qC)\Delta(0)-\frac{q(C+m)}3=0\;, \quad \Delta \text{ is even.}\]
Write
\[\omega_1 =(p+qC)\frac {\Delta'}{\sigma r} dx\wedge dy + \Delta'(q-pC)\frac{qe^1-pe^4}{p^2+q^2}\wedge dx +
\Delta e^{23}\;.\]
We already know that $\Delta(r)$ is even; hence,
\[\Delta'(r)=\frac{d\Delta}{dt}(r)=\mp\frac{d\Delta}{dr}(r)\]
is odd, implying that $\omega_1$ is smooth. On the other hand, requiring that $\omega_1$ be non-degenerate gives
\[\Delta''(0)\neq 0\;,\quad p+qC\neq 0\;.\]
Replacing $\xi$ with $-\xi$ if necessary, we can assume that $p+qC$ is positive.
The metric tensor can be decomposed  as
\begin{multline*} dx\otimes
dx+s(dy\otimes dy) + f\left(dy\odot\frac{qe^1-pe^4}{p^2+q^2}\right) +
g\left(\frac{qe^1-pe^4}{p^2+q^2}\right)^2 +\\+ (h^2+c^2)\left(e^2\otimes e^2\right)+
(b^2+k^2)\left(e^3\otimes e^3\right) + (hb+ck)e^2\odot e^3
\end{multline*}
whence
\[s=\frac1{\sigma^2 r^2}\left(\left(2(p+qC)\Delta-\frac{q(C+m)}3\right)^2+(p+qC)^2(\Delta')^2\right)\]
must be even and $s(0)=1$. The former condition is automatic, whereas the latter gives
\[\abs{\Delta''(0)}=\frac\sigma{p+qC}\;.\]
Moreover $f$ must be odd, $g$ must be even, and $g$, $h^2+k^2+c^2+b^2$ must be non-zero at the origin; all these conditions follow from those we have already obtained. The rest of the statement follows from the fact that $e^2\otimes e^2+e^3\otimes e^3$ is $H$-invariant, whereas
\[\langle e^2\odot e^3, e^2\otimes e^2-e^3\otimes e^3\rangle\]
is isomorphic to $V_{2p}$.
\end{proof}

We can now prove the main result of this section. We shall use the identification
\[T^2=\bigl\{(\exp 2\theta e_1, \exp \psi e_4)\in \SU(2)\times\LieG{U}(1)\bigr\}\cong\left\{(e^{i\theta}, e^{i\psi})\in\LieG{U}(1)\times\LieG{U}(1)\right\}\;.\]
\begin{theorem}
\label{thm:Compactification}
There is no solution of \eqref{eqn:EvolutionGivingNothing} that defines an Einstein-Sasaki metric on a compact manifold.
The Einstein-Sasaki structure \eqref{eqn:ExplicitMetric} extends to an invariant Einstein-Sasaki structure on the compact cohomogeneity one manifold with diagram
\eqref{eqn:Diagram}
if and only if:
\begin{itemize}
 \item  $A=0$,  $K$ is a finite subgroup of $T^2$ intersecting $(H_+)_0$ in a group of order $\sigma$,
\begin{gather*}
H_+=\left\{\left(e^{\frac12(\sigma+mq)i\theta},e^{qi\theta}\right)\in T^2,\theta\in\R\right\}\cdot K\;,\\
H_-=(\SU(2)\times\{1\})\cdot K\;, \end{gather*}
where $q$ and $\frac12(qm+\sigma)$ are coprime integers, and $\sigma = -\frac13(C+m)q$.
The resulting
Riemannian manifold is locally isometric to  $S^5$.
\item
$A<0$, $K$ is a finite subgroup of $T^2$ intersecting $(H_\pm)_0$ in a group of order $\sigma_\pm$;
\[H_\pm=\left\{\left(e^{\frac12(\sigma_\pm+mq_\pm) i\theta},e^{q_\pm i\theta}\right)\in T^2,\theta\in\R\right\}\cdot K,\]
where $q_\pm$ and $\frac12(q_\pm m+\sigma_\pm)$ are coprime integers,
and the equation
\begin{equation}
 \label{eqn:RootBloodyRoot}A+\Delta^2-4\Delta^3=0
\end{equation}
has two distinct roots $\Delta_-$, $\Delta_+$ such that
\begin{equation*}
 \frac{q_\pm}{\sigma_\pm}=\frac1{C+m}\,\frac{6\Delta_\pm}{1-6\Delta_\pm}\;.
\end{equation*}
Then the Einstein-Sasaki structure is quasi-regular if $C$ is rational, and irregular otherwise.
\end{itemize}
\end{theorem}
\begin{proof}
In \eqref{eqn:EvolutionGivingNothing}, we can assume that at some $t=t_0$
\[h,k,\Delta'>0,\quad b=0,\quad (k-h,c)\neq(0,0)\;.\]
Introduce the variables
\[u=h+k\;,\quad v=h-k\;,\quad z=b+c\;,\quad w=b-c\;;\]
by above, either $u,v\neq 0$ at $t=t_0$ or we can apply a transformation of type \eqref{eqn:U1Action} to obtain $u,v\neq 0$.
Then about $t_0$ \eqref{eqn:EvolutionGivingNothing} reads
\begin{equation}
\label{eqn:DeltaNothing}\Delta''+6\Delta=-1-\Delta' \frac{z'}z=1-\Delta' \frac{w'}w=1-\Delta' \frac{u'}u=-1-\Delta' \frac{v'}v
\end{equation}
Every solution of \eqref{eqn:DeltaNothing} satisfies $w=\lambda u$ and $z=\mu v$, where $\lambda$ and $\mu$
are constants. Thus
\[\Delta=\frac14(u^2-v^2-z^2+w^2)=\left(\frac{1+\lambda^2}4\right)u^2-\left(\frac{1+\mu^2}4\right)v^2\;.\]
This motivates us to set $U=\frac{\sqrt{1+\lambda^2}}2u$, $V=\frac{\sqrt{1+\mu^2}}2v$, transforming \eqref{eqn:DeltaNothing} into
\begin{equation}
\label{eqn:Delta}\Delta''+6\Delta=1-\Delta'\frac{U'}U=-1-\Delta'\frac{V'}V,\quad
\Delta'=\frac2{\frac{U'}U-\frac{V'}V},\quad \Delta=U^2-V^2\;.
\end{equation}
These equations only make sense at points where $U,V\neq 0$. However, $U$ cannot vanish in $(t_-,t_+)$ because $\Delta$ is positive, and if $V$ were zero at some $t_0$ in $(t_-,t_+)$ then it would be zero on all the interval, as one can see by applying a transformation of type \eqref{eqn:U1Action} and reducing to case (\emph{ii}) of Section~\ref{sec:Explicit}.

Suppose that $H$ is one-dimensional. Since $C=0$, Lemma~\ref{lemma:Extends} implies that $p>0$, and
\[hb+ck=\frac12(\lambda+\mu)uv\;,\quad  h^2+c^2-b^2-k^2=(1-\lambda\mu)uv\]
must vanish at the origin. So, $U(0)V(0)=0$. Since, again by Lemma~\ref{lemma:Extends}, $\Delta$ is a smooth
function with $\Delta(0)>0$, it follows that $V(0)=0$, $U(0)\neq 0$. Now observe that
\begin{equation}
\label{eqn:Squares}
h^2+k^2+b^2+c^2=4V^2+2\Delta=4U^2-2\Delta\;,
\end{equation}
is even by Lemma~\ref{lemma:Extends}, so that $U^2$ and $V^2$ are also even. This implies that $U$, $V$ are smooth functions, and $U'/U$ vanishes at $0$. By \eqref{eqn:Delta}, it follows that
\[\lim_{r\to 0}\frac{V'}V=\lim_{r\to 0}-\left(\frac{U'}U-\frac{V'}V\right)=\lim_{r\to 0}-\frac2{\Delta'}=-\infty\;.\]
On the other hand, since $V(r)$ is non-zero for all $r>0$ we have
\begin{equation*}
 \lim_{r\to 0} \frac{\frac{dV}{dr}}{V}\geq 0\;.
\end{equation*}
Since $dr=\mp dt$, we conclude that this possibility may only occur at $t_+$.

Suppose now that $H$ has dimension three. Then  Lemma~\ref{lemma:Extends} and \eqref{eqn:Squares} imply that
$V^2/r^2$, $U^2/r^2$ are  even. Dividing \eqref{eqn:Squares} by $r^2$ and evaluating at the origin, we find
$U^2/r^2=\frac14$, $V^2/r^2=0$ at the origin.  It follows from
\eqref{eqn:Delta} and Lemma~\ref{lemma:Extends} that
\[\frac12=-1-\lim_{r\to 0}\Delta'\frac{V'}V=-1\pm\frac12 \lim_{r\to 0} \frac{rV'}V\]
and so we reach the contradiction
\[0\leq \lim_{r\to 0}\frac{r\frac{dV}{dr}}{V}=\mp \lim_{r\to 0}\frac{rV'}{V}=-3\;.\] This completes the
proof of the first part of the theorem.

\smallskip
Now consider a solution of \eqref{eqn:ExplicitMetric}. By Lemma~\ref{lemma:Extends}, at the origin either
$h=0$ or \mbox{$(h^2)'=0$}. The first possibility can only occur if $A=0$. In this case,
\eqref{eqn:ChangeOfVariable} gives
\[h'=\sqrt{\frac14-h^2}\;,\]
and thus $\Delta=h^2$ is a smooth even function such that  $\Delta''(0)=\frac12$. Then, by
Lemma~\ref{lemma:Extends}, $H_0=\SU(2)\times\{1\}$ and the structure extends across $G/H$. Notice that one
cannot have $h=0$ at both special orbits, because otherwise $h'$ would have to vanish at some point in
between. On the other hand, $h'$ vanishes at both special orbits if and only if $A$ is negative.

Suppose that $h'=0$ at the origin. It follows from \eqref{eqn:ChangeOfVariable} that $\Delta$ is
automatically even, and Lemma~\ref{lemma:Extends} only requires that, at the origin,
\[\Delta=\frac{q(C+m)}{6(p+qC)}, \quad\frac\sigma{p+qC}= \abs{1-6\Delta}\;,\]
where we have used a formula analogous to \eqref{eqn:DeltaNothing} to express $\Delta''$ in terms of
$\Delta$. Since $\sigma$ is determined up to sign, we can write
\begin{equation}
\label{eqn:PAndQ} q=\frac{\sigma}{C+m}\frac{6\Delta(0)}{1-6\Delta(0)}\;,\quad p=qm+\sigma\;.
\end{equation}
The condition $h'=0$ means that $\Delta(0)$ is a root of
\eqref{eqn:RootBloodyRoot}.

Thus, if $A=0$, $\Delta(t_+)=\frac14$ and \eqref{eqn:PAndQ} gives
\[q=-\frac{3\sigma}{C+m}\;;\]
the remaining conditions on $p$, $q$ and $\sigma$ follow from Lemma~\ref{lemma:Extends}.
A straightforward computation shows that the
Riemannian metric has constant sectional curvature, and therefore the Riemannian manifold we obtain is
locally isometric to $S^5$.

If $A$ is non-zero, both special orbits are three-dimensional and satisfy
\eqref{eqn:PAndQ}, which together with Lemma~\ref{lemma:Extends} determines the required conditions. Finally, the statement about regularity follows from the fact that the characteristic vector field is $Ce_1-e_4$, up to an invariant function.
\end{proof}
It was proved  in \cite{GauntlettMartelliSparksWaldram} that there are countably infinite values of
$A$ in the interval \mbox{$-\frac1{108}<A<0$} for which \eqref{eqn:RootBloodyRoot} has roots $\Delta_+$,
$\Delta_-$ with $\Delta_+-\Delta_-\in\Q$. One can show, using \eqref{eqn:RootBloodyRoot}, that for any such value of $A$ and any choice of $m$ there exists a
value of $C$ satisfying the hypotheses of Theorem~\ref{thm:Compactification}. Moreover, there also exist
infinite values of $A$ for which both $\Delta_+$ and $\Delta_-$ are rational, so that $C$ is also rational
and the metric is quasi-regular.

\begin{remark} It follows from \eqref{eqn:ExplicitMetric} that the right action of $T^2$ on $G$ preserves the
$\LieG{U}(2)$-structure. One can then replace $G$ with a subgroup
\[\SU(2)\times\LieG{U}(1)\subset G\times T^2\]
that acts transitively on $G$ preserving the Einstein-Sasaki $\SU(2)$-structure.  Since \mbox{$T^2\supset K,H_\pm$} is abelian, this action is well defined on the cohomogeneity one manifold \eqref{eqn:Diagram}. Thus, we
discover that the cohomogeneity one action can be assumed to preserve the $\SU(2)$-structure, which amounts
to setting $m=0$ in \eqref{eqn:ExplicitMetric}.
\end{remark}
By the above remark, there is no loss of generality in assuming $m=0$.  Thus, we recover the Einstein-Sasaki
metrics of \cite{GauntlettMartelliSparksWaldram}. Indeed, introduce coordinates $(\theta,\psi,\phi)$ on
$\SU(2)=S^3$ as in \eqref{eqn:Euler}, a coordinate $\alpha$ on $\LieG{U}(1)$ such that $e^4=d\alpha$, and a
coordinate $y$ on $(t_-,t_+)$ such that $y=1-6\Delta$. Set
\[\beta=-\psi+C\,\alpha\;,  \quad wq=\left(\frac{2(108 A+1-3y^2+2y^3)}{1-y}\right)\;.\]
Then our  family of local Einstein-Sasaki metrics can be written as
\begin{multline*}
\frac{1-y}6 (d\theta^2 +\sin^2\theta d\phi^2)+\frac{1}{wq} dy^2+\frac1{36} wq(d\beta+\cos\theta d\phi)^2 +\\
+\frac19 \left(d\psi-\cos\theta d\phi + y(d\beta+\cos\theta d\phi)\right)^2
\end{multline*}
which is the local form of a $Y^{p,q}$ metric if $A<0$, and the standard metric on the sphere if $A=0$.

\section{Classification}
We can finally prove that the only cohomogeneity one Einstein-Sasaki $5$-manifolds are the $Y^{p,q}$.
\begin{theorem}
\label{thm:Classification}
Let $\ev{M}$ be a compact, simply-connected Einstein-Sasaki $5$-manifold on which the group of isometries acts with cohomogeneity one. Then $\ev{M}$ can be represented by the diagram \eqref{eqn:Diagram}, where
\begin{align*}
K&=\left\langle \left(-1,e^{2\pi iq_+/\sigma_+}\right),\left(-1,e^{2\pi iq_-/\sigma_-}\right)\right\rangle\subset T^2,\\
H_\pm&=\left\{\left(e^{\frac12\sigma_\pm i\theta},e^{q_\pm i\theta}\right)\in T^2,\theta\in\R\right\}\cdot K,
\end{align*}
the integers $q_\pm$, $\sigma_\pm$ satisfy  Theorem~\ref{thm:Compactification}, and the metric is  given by \eqref{eqn:ExplicitMetric}. Moreover, $\ev{M}$ is diffeomorphic to $S^2\times S^3$.
\end{theorem}
\begin{proof}
Since the group of isometries of $S^5$ acts in a homogeneous way, $\ev{M}$ is not isometric to $S^5$. Thus, by Theorem~\ref{thm:Compactification} and the subsequent remark $\ev{M}$ can be represented by the diagram $\tilde K\subset{\tilde H_\pm}\subset G$, where
$\tilde K$ is a finite subgroup of $T^2$ intersecting $(\tilde H_\pm)_0$ in a group of order $\sigma_\pm$, and
\[\tilde H_\pm=\left\{\left(e^{\frac12\sigma_\pm i\theta},e^{q_\pm i\theta}\right)\in T^2,\theta\in\R\right\}\cdot \tilde K\;.\]
Now define $K$, $H_\pm$ as in the statement; then $\tilde H_\pm/H_\pm\cong \tilde K/K$. Thus, up to an equivariant covering map $\ev{M}$ is represented by the diagram $K\subset{H_\pm}\subset G$. Since $\ev{M}$ is simply connected, this covering map is actually a diffeomorphism.

Recall from Section~\ref{sec:CohomogeneityOne} that the diagram determines $\ev{M}$ only up to a glueing map. However, in our case $K$ is contained in the center of $G$, and so  the diagram does determine $\ev{M}$ up to equivariant diffeomorphism.

To determine the topology of $\ev{M}$, we make use of Smale's classification theorem~\cite{Smale}; indeed, the existence of an Einstein-Sasaki metric implies that $\ev{M}$ is spin,  whence it suffices to
prove that
\begin{equation*}
 H_2(\ev{M})=\Z,\quad \pi_1(\ev{M})=0\;.
\end{equation*}
Let $D_\pm$ be the
$2$-dimensional disc, and write $\ev{M}=U_+\cup U_-$, where $U_\pm=G\times_{H_\pm}  D_\pm$. Then $G/H_\pm$
is a retract of $U_\pm$, whereas $G/K$ is a retract of $U_+\cap U_-$. Define loops
$\gamma_4,\gamma_\pm\colon S^1\to G$ by
\[\gamma_4(e^{it})=\exp(te_4),\quad  \gamma_\pm(e^{it})=\exp(t\sigma_\pm e_1+tq_\pm e_4)\;.\]
It is not difficult to show that
\begin{gather*}
\pi_1(G/K)=\frac{\Z\gamma_+\oplus\Z\gamma_-\oplus\Z\gamma_4}{\langle\sigma_+\gamma_+ -q_+\gamma_4,\sigma_-\gamma_--q_-\gamma_4\rangle}\;,\\
\pi_1(G/H_\pm)=\frac{\Z\gamma_-\oplus\Z\gamma_4}{\langle \sigma_\mp\gamma_\mp- q_\mp\gamma_4,
q_\pm\gamma_4\rangle}\;,\quad H_2(G/H_\pm)=0 \;.
\end{gather*}
 By the Van-Kampen theorem, it follows that
\[\pi_1(\ev{M})=\frac{\Z\gamma_4}{\langle q_+\gamma_4,q_-\gamma_4\rangle}\;,\]
and so $\ev{M}$ is simply connected if and only if  $q_+$ and $q_-$ are coprime. On the other hand the
Mayer-Vietoris sequence gives
\[
\xymatrix{
 H_2(U_+)\oplus H_2(U_-)\ar[r] & H_2(M) \ar[r]&H_1(U_+\cap U_-)\ar^j[r] & H_1(U_+)\oplus H_1(U_-)
}\] By above, we conclude \[H_2(\ev{M})\cong\ker j=\Z(q_+q_-\gamma_4)\cong\Z\;,\] where we have used the fact that
$q_+$ and $q_-$ are coprime.
\end{proof}

\noindent\textbf{Acknowledgements.}
 I would like to thank A. Ghigi and S. Salamon for helpful discussions. I am also in debt to K.~Galicki and the referee for some useful comments on  earlier versions of this paper.
\bibliographystyle{plain}
\bibliography{cohomo}
\end{document}